\documentclass{amsart}

\usepackage{amssymb,amsmath,amsfonts}
\usepackage[all]{xy}

\title[Generating series for Hilbert schemes of points]{Generating series in the cohomology of Hilbert schemes of points on
surfaces}
\author{Samuel Boissi\`ere \and Marc A.~Nieper-Wi\ss kirchen}

\date{\today}
\address{Samuel Boissi{\`e}re, Laboratoire J.A.Dieudonn\'e UMR CNRS 6621,
         Universit\'e de Nice Sophia-Antipolis, Parc Valrose, 06108 Nice, France}
\email{sb@math.unice.fr}
\address{Marc A.~Nieper-Wi\ss kirchen, Institut f\"ur Mathematik\\Johannes-Gutenberg Universit\"at\\
 55099 Mainz\\Germany}
\email{nieper@mathematik.uni-mainz.de}

\newcommand{\IC}{\mathbb C}
\newcommand{\IH}{\mathbb H}
\newcommand{\IQ}{\mathbb Q}

\newcommand{\vac}{|0\rangle}
\newcommand{\unit}{|1\rangle}

\newcommand{\cO}{\mathcal O}
\newcommand{\cP}{\mathcal P}

\newcommand{\ie}{{\it i.e. }}

\newcommand{\kc}{\mathfrak{c}}
\newcommand{\kch}{\mathfrak{ch}}
\newcommand{\kd}{\mathfrak{d}}
\newcommand{\kf}{\mathfrak{f}}
\newcommand{\kg}{\mathfrak{g}}

\newcommand{\kq}{\mathfrak{q}}

\newcommand{\kG}{\mathfrak{G}}
\newcommand{\kS}{\mathfrak{S}}

\DeclareMathOperator{\Hilb}{Hilb}
\DeclareMathOperator{\Ch}{Ch}
\DeclareMathOperator{\ch}{ch}
\DeclareMathOperator{\CC}{C}
\DeclareMathOperator{\cc}{c}
\DeclareMathOperator{\Se}{S}
\DeclareMathOperator{\Td}{Td}
\DeclareMathOperator{\td}{td}
\DeclareMathOperator{\ad}{ad}
\DeclareMathOperator{\rk}{rk}
\DeclareMathOperator{\End}{End}

\newtheorem{theorem}{Theorem}[section]
\newtheorem{formula}[theorem]{Formula}
\newtheorem{proposition}[theorem]{Proposition}
\newtheorem{remark}[theorem]{Remark}

\begin{document}

\begin{abstract}
In the study of the rational cohomology of Hilbert schemes of points on a
smooth surface, it is particularly interesting to understand the characteristic
classes of the tautological bundles and the tangent bundle. In this note we
pursue this study. We first collect all results appearing separately in the
literature and prove some new formulas using T. Ohmoto's results on orbifold
Chern classes on Hilbert schemes. We also explain the algorithmic counterpart
of the topic: The cohomology space is governed by a vertex algebra that can be
used to compute characteristic classes. We present an implementation of the
vertex operators in the rewriting logic system {\sc Maude} and address
observations and conjectures obtained after symbolic computations.
\end{abstract}

\maketitle

\tableofcontents

\part{Preliminaries}
\label{p:prel}

\section{Introduction}

Let $S$ be a smooth quasi-projective complex surface and $n\geq 0$ an integer.
The \emph{Hilbert scheme} of $n$ points on $S$, denoted by $S^{[n]}$, is the
moduli space of generalized $n$-tuples on $S$, \ie zero-dimensional subschemes
of length $n$ on $S$. $S^{[n]}$ is smooth of complex dimension $2n$. When
working in the \emph{rational} cohomology of Hilbert schemes, it is usual to
consider the total Hilbert scheme:
\begin{displaymath}
\Hilb(S):=\coprod_{n\geq 0} S^{[n]}
\end{displaymath}
whose total cohomology space is:
\begin{displaymath}
\IH_S:=\prod_{n\geq 0}\bigoplus_{i=0}^{4n}H^i(S^{[n]}).
\end{displaymath}
The space $\IH_S$ is completed bigraded by \emph{conformal weight} $n$ and
\emph{cohomological degree} $i$.

Consider the tangent bundle $T_S^n$ on $S^{[n]}$ of rank $2n$. Any
characteristic class $\phi$ (such as the Chern class, the Segre class, the Todd
class, the Chern character) can be applied to this bundle and gives, after
summing up over all values of $n$, an element:
\begin{displaymath}
\Phi(S):=\sum_{n\geq 0} \phi\left(T_S^n\right) \in \IH_S.
\end{displaymath}
There exists a universal way of presenting these characteristic classes,
involving the canonical class $K_S$ and the Euler class $e_S$ of $S$, depending
on universal constants independent of $S$, indexed by partitions of integers.

\bigskip

Similarly, any vector bundle $F$ of rank $r$ on $S$ defines in a tautological
way a vector bundle $F^{[n]}$ of rank $n\cdot r$ on $S^{[n]}$, called a
\emph{tautological bundle}. Any characteristic class $\phi$ can be applied to
these bundles and gives an element:
\begin{displaymath}
\Phi(F):=\sum_{n\geq 0} \phi\left(F^{[n]}\right) \in \IH_S.
\end{displaymath}
As before, there exists a universal way of presenting these characteristic
classes, involving the corresponding characteristic class $\phi(F)$, the
canonical class $K_S$ and the Euler class $e_S$ of $S$, depending on universal
coefficients independent of $F$ and $S$, indexed by partitions of integers.

\bigskip

In this paper, we are interested in the effective computation of these
universal coefficients in both situations. Since the constants are independent
of the surface, one may evaluate each formula on well-chosen surfaces and use
appropriate tools to compute new coefficients step by step. This is one of the
main tricks used to obtain the values of some series of coefficients, together
with manipulations of vertex algebra operators. When these methods fail, one
can make use of a suitable computer program to get information on the missing
values.

In Section 2 we recall some basics on the vertex algebra structure of the total
cohomology space $\IH_S$: natural vertex operators and commutation relations
between them lead to effective algorithms to compute classes in $\IH_S$. Part
\ref{p:taut} is devoted to the case of the tautological bundles and Part
\ref{p:tangent} to the tangent bundle: We explain the general shape of the
formulas, then present the state of achievement in the seek of the universal
constants and prove some new results, and finally we produce new values
obtained by an implementation of the vertex algebra structure with {\sc Maude}
\cite{Maude}: the vertex operators are directly defined with their commutation
rules, making the program clear to understand from the mathematical point of
view (Part \ref{p:maude}). To our knowledge, {\sc Maude} has not been used for
huge algebraic computations and simplifications so far. The nature of {\sc
Maude} as a rewriting system will make it easy to prove the correctness of the
implemented algorithms.

\section{Notation}

We recall here some notation and classical constructions. For further details,
refer to \cite{B,BNW,Lehn,LehnLN}.

\subsection{Combinatorics}\text{}

Let $\cP$ be the set of all partitions of integers. A partition
$\lambda=(\lambda_1,\lambda_2,\cdots)$ is likewise denoted by
$\lambda=(1^{m_1},2^{m_2},\ldots)$, and we define:
\begin{align*}
\ell(\lambda)&:=\sum_{i\geq 1} m_i \quad\text{(the length)}\\
|\lambda|&:=\sum_{i\geq 1} i m_i \quad\text{(the weight)}\\
||\lambda||&:=\sum_{i\geq 1} i^2 m_i \\
\lambda!&:=\prod_{i\geq 1}m_i!
\end{align*}
We often write $n$ instead of $(n)=(n^1)$.

\subsection{Vertex operators}\text{}

A linear operator $\kf\in \End(\IH_S)$ is \emph{homogeneous} of conformal
weight $u$ and cohomological degree $v$ if for any $n$ one has
$\kf\left(H^i(S^{[n]})\right)\subset~ H^{i+v}(S^{[n+u]})$. The
\emph{commutator} of two homogeneous operators $\kf,\kg$ is defined by:
\begin{displaymath}
[\kf,\kg]:=\kf\circ\kg-(-1)^{|\kf|\cdot|\kg|}\kg\circ \kf
\end{displaymath}
where $|\cdot|$ denotes the cohomological degree.

The total cohomology space $\IH_S$ is computed with the help of
\emph{Nakajima's creation operators} \cite{Nak}:
\begin{displaymath}
\kq_k:H^*(S)\longrightarrow \End(\IH_S),\quad k\geq 1.
\end{displaymath}
For $\alpha\in H^i(S)$, the operator $\kq_k(\alpha)$ has conformal weight $k$
and cohomological degree $2(k-1)+i$. We shall make use of the following
abbreviation: For $k\geq 1$, the push-forward induced by the diagonal inclusion
$\Delta^k:S\rightarrow S^k$ gives a map $\Delta^k_!:~H^*(S)\rightarrow~
H^*(S)\otimes~\cdots\otimes~ H^*(S)$. For $\Delta^k_!\alpha=~\sum_i
\alpha_{i,1}\otimes~\cdots\otimes~\alpha_{i,k}$ and
$\lambda=(\lambda_1,\ldots,\lambda_k)$ a partition of length $k$ we set (see
\cite{B,BNW}):
\begin{displaymath}
\kq_\lambda(\alpha):=(\kq_{\lambda_1}\circ\cdots\circ\kq_{\lambda_k})\Delta^k_!(\alpha):=\sum_i\kq_{\lambda_1}(\alpha_{i,1})\circ
\cdots\circ\kq_{\lambda_k}(\alpha_{i,k}).
\end{displaymath}

The unit in $H^0(S^{[0]})\cong \IQ$ is denoted $\vac$ and is called the
\emph{vacuum} of $\IH_S$. Evaluations on the vacuum, like
$\kq_\lambda(\alpha)\vac$, provide natural classes in $\IH_S$. The unit in
$\IH_S$ for the cup product is given by:
$$
\unit:=\sum_{n\geq 0}1_{S^{(n]}}=\mathrm{e}^{\kq_1(1_S)}\vac.
$$

\subsection{Tautological classes}\text{}

Let $\Xi^n_S$ be the universal family on $S^{[n]}$ with the projections:
\begin{displaymath}
\xymatrix{\Xi^n_S\ar@{^(->}[r]&S^{[n]}\times S\ar[d]_p \ar[r]^-q& S \\&
S^{[n]}}
\end{displaymath}
Let $F$ be a locally free sheaf on $S$. For any $n\geq 0$, the associated
\emph{tautological bundle} on $S^{[n]}$ is defined by:
\begin{displaymath}
F^{[n]}:=p_*\left(\cO_{\Xi^n_S}\otimes q^*F\right).
\end{displaymath}
It has rank $n\cdot \rk(F)$. This construction extends naturally to a
well-defined group homomorphism:
\begin{displaymath}
-^{[n]}:K(S)\rightarrow K(S^{[n]})
\end{displaymath}
where $K(\cdot)$ denotes the \emph{rational} Grothendieck group generated by
locally free sheaves. For $u\in K(S)$, let $\kc(u)$ and $\kch(u)$ be the linear
operators acting for any $n\geq 0$ on $H^*(S^{[n]})$ by the cup-product
multiplication by the total Chern class $\cc(u^{[n]})$ and the total Chern
character $\ch(u^{[n]})$ respectively.

In particular, taking $F=\cO_S$ one defines a linear operator $\kd\in \End(\IH^S)$ by:
\begin{displaymath}
\kd(x):=c_1(\cO_S^{[n]})\cdot x,\quad \forall x\in H^*(S^{[n]}).
\end{displaymath}
The \emph{derivative} of a linear operator $\kf\in \End(\IH_S)$ is defined by
$\kf':=[\kd,\kf]$.

By analogy with the construction of tautological bundles, one defines a linear
operation $-^{[n]}:H^*(S)\rightarrow H^*(S^{[n]})$: For any cohomology class
$\gamma\in H^*(S)$ we set
\begin{displaymath}
\gamma^{[n]}:=p_*\left(\ch(\cO_{\Xi^n_S})\cdot q^*\td(S)\cdot q^*\gamma\right)
\end{displaymath}
where $\td(S)$ denotes the Todd class of the tangent bundle on $S$ and we
define an operator $\kG(\gamma)\in \End(\IH^S)$ acting on $H^*(S^{[n]})$ by
multiplication by $\gamma^{[n]}$.

Similarly, denoting by $-^\vee:K(S)\rightarrow K(S)$ the natural ring
involution taking the dual of a vector bundle, for $\gamma\in H^*(S)$ we set
\begin{displaymath}
\left(\gamma^{[n]}\right)^\vee:=p_*\left(\ch(\cO_{\Xi^n_S}^\vee)\cdot
q^*\td(S)\cdot q^*\gamma\right)
\end{displaymath}
and define an operator $\kG^\vee(\gamma)\in \End(\IH^S)$ acting on
$H^*(S^{[n]})$ by multiplication by $\left(\gamma^{[n]}\right)^\vee$.

\subsection{Cohomology of $S$}\text{}

In the rational cohomology ring of $S$, we denote the unit by $1_S\in H^0(S)$,
the canonical class (\ie the first Chern class of the cotangent bundle on $S$)
by $K_S\in~ H^2(S)$ and the Euler class of $S$ (\ie the second Chern class of
the tangent bundle on $S$) by $e_S\in H^4(S)$.

\part{Characteristic classes of tautological bundles}
\label{p:taut}

In this first part, we study the characteristic classes of a tautological
bundle on $S^{[n]}$ obtained from a vector bundle $F$ on $S$. We first recall
the theoretic results.

\section{Shape of the formulas}

For tautological bundles, the best general result is obtained for the Chern
character:

\begin{theorem}[Boissi\`ere {\cite{B}}, Boissi\`ere \& Nieper-Wi\ss kirchen {\cite{BNW}}]\label{ChernCharTaut} There are unique
rational constants $\alpha_\lambda,\beta_\lambda,\gamma_\lambda,\delta_\lambda$
such that for each surface $S$ and each vector bundle $F$ on $S$, the
generating series of the Chern characters of the tautological bundles of $F$ is
given by:
\begin{align*}
\Ch(F):=\kch(F)\unit
=\Big(\sum_{\lambda\in \cP}&\alpha_\lambda\kq_\lambda(\ch(F))+\beta_\lambda\kq_\lambda(e_S\ch(F))\\
&+\gamma_\lambda\kq_\lambda(K_S\ch(F))+\delta_\lambda\kq_\lambda(K_S^2\ch(F))\Big)\unit.
\end{align*}
\end{theorem}

For multiplicative characteristic classes, the general shape is similar. Let
$\phi$ be a multiplicative characteristic class. Let $B$ the polynomial ring
$\IQ[r,c_1,c_2,K,e]$ truncated from degree $5$ onwards where $\deg(r)=0$,
$\deg(c_1)=2$, $\deg(c_2)=4$, $\deg(K)=2$, $\deg(e)=4$.

\begin{theorem}[Boissi\`ere {\cite{B}}, Boissi\`ere \& Nieper-Wi\ss kirchen
\cite{BNW}] There are unique elements $u_\lambda^\phi\in B$ such that for each
surface $S$ and each vector bundle $F$ on $S$, the generating series of the
$\phi$-classes of tautological bundles is given by:
\begin{displaymath}
\Phi(F):=\sum_{n\geq 0}
\phi(F^{(n]})=\exp\left(\kq_\lambda(u_\lambda^\phi(\rk(F),c_1(F),c_2(F),K_S,e_S)\right)\vac.
\end{displaymath}
\end{theorem}

As we shall see in Section \ref{ss:rank1}, the formula for the total Chern
class simplifies considerably when specialized at $r=1$. As an example of the
complexity in higher rank, we shall give in Section \ref{ss:rankhigh} the first
terms of the linear series for $r\geq 2$.

\section{The Chern Character}

In the determination of the constants in Theorem \ref{ChernCharTaut}, a lot of
information has already been obtained, and the result is complete for surfaces
with trivial canonical class.

\subsection{The $(1_S)$- and $(e_S)$-series}\text{}

All coefficients $\alpha_\lambda$ and $\beta_\lambda$ are known:

\begin{formula}[Li, Qin \& Wang {\cite[Corollary 4.8]{LQW}}] It is:
\begin{align*}
\alpha_\lambda &= \frac{(-1)^{|\lambda|-1}}{\lambda!\cdot |\lambda|!},\\
\beta_\lambda &=
\frac{(-1)^{|\lambda|}}{\lambda!|\lambda|!}\cdot\frac{|\lambda|+||\lambda||-2}{24}.
\end{align*}
\end{formula}

\subsection{The $(K_S)$- and $(K_S^2)$-series}\text{}
\label{chtautimplement}

For the series concerning the canonical class, Li, Qin \& Wang \cite[Corollary
4.8]{LQW} write the still unknown constants $\gamma_\lambda$ and
$\delta_\lambda$ as
\begin{displaymath}
\gamma_\lambda = \frac{(-1)^{|\lambda|}}{\lambda!\cdot |\lambda|!}\cdot
g_{1}(\lambda\cup 1^{|\lambda|})\qquad \text{and}\qquad \delta_\lambda =
\frac{(-1)^{|\lambda|}}{\lambda!|\lambda|!}\cdot g_{2}(\lambda\cup
1^{|\lambda|}).
\end{displaymath}
where the functions $g_1,g_2$ only depend on the partition and $\lambda\cup
1^{|\lambda|}$ means that one adds $|\lambda|$ to the multiplicity of $1$ in
$\lambda$.

In order to get information about these functions we implement the commutation
relation of Lehn \cite[Theorem 4.2]{Lehn}:
\begin{displaymath}
[\kch(F),\kq_1(1_S)]=\exp(\ad \kd)(\kq_1(\ch(F)).
\end{displaymath}
This gives the following recursive formula (see \cite[\S 3.3]{B}):
\begin{displaymath}
\ch\left(F^{[n]}\right)=\frac{1}{n}\kq_1(1_S)\ch\left(F^{[n-1]}\right)+
\sum\limits_{\nu=0}^{2n}\frac{\kq^{(\nu)}_1(\ch(F))}{\nu!}\frac{\kq_1(1_S)^{n-1}}{(n-1)!}
\vac.
\end{displaymath}
Computations with {\sc Maude} \cite{Maude} (see Part \ref{p:maude}) give the
following values (for each value, we extract the factor $
\frac{1}{\lambda!\cdot |\lambda|!}$):

\begin{displaymath}
\begin{array}{|c||c|c|c|c|c|c|}
\hline \lambda & (1) & (1,1) & (2) & (1,1,1) & (2,1) & (3) \\
\hline \gamma & 0 & -\frac{1}{4}\cdot\frac{1}{3} & -\frac{1}{2}\cdot 1 &
\frac{1}{36}\cdot\frac{3}{5} & \frac{1}{6}\cdot\frac{3}{4} & \frac{1}{6}\cdot 1 \\
\hline \delta & 0 & -\frac{1}{4}\cdot\frac{1}{12} &
-\frac{1}{2}\cdot\frac{1}{6} & \frac{1}{36}\cdot\frac{7}{30} &
\frac{1}{6}\cdot\frac{7}{20} & \frac{1}{6}\cdot\frac{7}{12}
\\\hline
\end{array}
\end{displaymath}

\begin{displaymath}
\begin{array}{|c||c|c|c|c|c|}
\hline \lambda & (1,1,1,1) & (2,1,1) & (2,2) & (3,1) & (4) \\
\hline \gamma & - \frac{1}{576}\cdot\frac{29}{35} & -\frac{1}{48}\cdot\frac{29}{30} & -\frac{1}{48}\cdot\frac{11}{10} & -\frac{1}{24}\cdot\frac{6}{5} & -\frac{1}{24}\cdot\frac{3}{2} \\
\hline \delta & - \frac{1}{576}\cdot\frac{59}{140} &
-\frac{1}{48}\cdot\frac{59}{105} & -\frac{1}{48}\cdot\frac{43}{60} &
-\frac{1}{24}\cdot\frac{5}{6} & -\frac{1}{24}\cdot\frac{5}{4}\\\hline
\end{array}
\end{displaymath}

\begin{remark}
This computation shows in particular that the functions $g_1$ and $g_2$
are not integer-valued and one may suppose that they are always negative.
\end{remark}

\section{The Chern class}

\subsection{The rank 1 case}\text{}\label{ss:rank1}

In the case of a vector bundle $F$ of rank $1$, there is a complete answer to
the question of the determination of the universal constants for the Chern
class:

\begin{formula}[Lehn {\cite[Theorem 4.6]{Lehn}}] Let $L$ be a line bundle on
$S$. The generating series of the  Chern classes of tautological bundles take
the form:
\begin{displaymath}
\CC(L):=\kc(L)\unit=\exp \left(\sum_{k\geq 1}
\frac{(-1)^{k-1}}{k}\kq_{k}(\cc(L))\right)\vac.
\end{displaymath}
\end{formula}

In particular, for $\lambda=(k)$,
$u_\lambda^{\cc}=\frac{(-1)^{k-1}}{k}\cc(L)+(r-1)(\cdots)$ and
$u_\lambda^{\cc}=(r-1)(\cdots)$ otherwise.

In this formula, one sees that the only operators $\kq_\lambda$ occurring have
partitions with one part and that the invariants $K_S$ and $e_S$ do not appear.
It is not expected that this remains true for other characteristic classes nor
in higher rank. There are no similar formulas known for other multiplicative
characteristic classes.

\subsection{Bundles of higher rank}\text{}\label{ss:rankhigh}

For a bundle $F$ of rank $r\geq 1$, information on the beginning of the
universal formula of the Chern class is contained in the following result
obtained by specialization to the affine plane:

\begin{formula}[Boissi\`ere \& Nieper-Wi\ss kirchen {\cite[Theorem 4]{BNW}}]
\label{cherntauttrivial} Let $F$ be the trivial bundle of rank $r$ on $\IC^2$.
Then:
\begin{displaymath}
\CC(F)=\exp \left(\sum_{k\geq 1} \frac{(-1)^{k-1}}{k^2}\binom{r\cdot
k}{k-1}\kq_{k}(1_{\IC^2})\right)\vac.
\end{displaymath}
\end{formula}

In particular, this gives the complete answer for the degree zero part of the
$u_{(k)}$.

In order to get more information, we implement the commutation relation of Lehn
\cite[Theorem 4.2]{Lehn}:
\begin{displaymath}
\kc(F)\circ\kq_1(1_S)\circ\kc(F)^{-1}=\sum_{\nu,k\geq
0}\binom{\rk(F)-k}{\nu}\kq_1^{(\nu)}(c_k(F)),
\end{displaymath}
which gives the following recursion formula:
\begin{displaymath}
\cc(F^{[n]})=\frac{1}{n}\left(\sum_{\substack{0\leq k\leq \rk(F)\\0\leq \nu\leq \rk(F)-k}}\binom{\rk(F)-k}{\nu}\kq_1^{(\nu)}(c_k(F))\right)\cc(F^{[n-1]}).
\end{displaymath}
Computations with {\sc Maude} \cite{Maude} give the series inside the
exponential. For example, in the rank $2$ case the first $u_\lambda^{\cc}$ are
(computed in the basis $1,c_1,K,c_2,c_1^2,c_1K,K^2,e$):
\begin{displaymath}
\begin{array}{|c|c|}
\hline \lambda & u_\lambda^{\cc} \\\hline \hline (1) & 1+c_1+c_2 \\\hline
\hline (1,1) &  -\frac{1}{2}(1+c_1+c_2)\\
\hline (2) &  -\left(1+\frac{1}{2}(3c_1+K+2c_2+c_1^2+c_1K)\right)\\\hline
\hline (1,1,1) & \frac{1}{3}(1+c_1+c_2) \\
\hline (2,1) & 2+3c_1+K+2c_2+c_1^2+c_1K \\
\hline (3) &
\frac{5}{3}+\frac{10}{3}c_1+2K+2c_2+2c_1^2+3c_1K+\frac{2}{3}K^2-\frac{1}{3}e
\\\hline
\hline (1,1,1,1) & -\frac{1}{4}(1+c_1+c_2) \\
\hline (2,1,1) & -(3+\frac{1}{2}(9c_1+3K+6c_2+3c_1^2+3c_1K))\\
\hline (2,2) & -\frac{1}{4}(9+18c_1+10K+10c_2+11c_1^2+15c_1K+3K^2-e)\\
\hline (3,1) & -(5+10c_1+6K+6c_2+6c_1^2+9c_1K+2K^2-e)\\
\hline (4) & -\frac{1}{4}(14+35c_1+29K+20c_2+30c_1^2+58c_1K+22K^2-10e)\\\hline
\end{array}
\end{displaymath}

\begin{remark}\text{}
In order to recover Formula \ref{cherntauttrivial} for $S=\IC^2$, set $c_1=0$,
$c_2=0$, $K=0$, $e=0$ and keep only the partitions with one part (since the
diagonal push-forward is trivial in this case).
\end{remark}

We can simplify the combinatorial difficulties by making the following
assumptions: $S$ is an abelian surface and $F$ is a trivial bundle of rank $r$
over $S$. Concretely, in our computer program (Part \ref{p:maude}) we set
$c_1=0$, $c_2=0$, $e=0$ and $K=0$. Thus the elements $u_\lambda^{\cc}$ inside
the exponential restrict to their degree zero term $(u_\lambda^{\cc})_0$ for
each partition. We get the following results:
\begin{itemize}
\item Rank $2$:
\begin{displaymath}
\begin{array}{|c||c|c|c|c|c|c|}
\hline \lambda & (1) & (1,1) & (2) & (1,1,1) & (2,1) & (3)\\
\hline (u_\lambda^{\cc})_0 & 1 & -\frac{1}{2} & -1 & \frac{1}{3} & 2 &
\frac{5}{3}\\\hline
\end{array}
\end{displaymath}

\begin{displaymath}
\begin{array}{|c||c|c|c|c|c|}
\hline \lambda & (1,1,1,1) & (2,1,1) & (2,2) & (3,1) & (4) \\
\hline (u_\lambda^{\cc})_0 & -\frac{1}{4} & -3 & -\frac{9}{4} & -5 &
-\frac{7}{2} \\\hline
\end{array}
\end{displaymath}

\begin{displaymath}
\begin{array}{|c||c|c|c|c|c|c|c|}
\hline \lambda & (1,1,1,1,1) & (2,1,1,1) & (2,2,1) & (3,1,1) & (3,2) & (4,1) &
(5) \\
\hline (u_\lambda^{\cc})_0 & \frac{1}{5} & 4 & 9 & \frac{61}{6} & 12 & 14 &
\frac{42}{5}\\\hline
\end{array}
\end{displaymath}

\begin{align*}
&\begin{array}{|c||c|c|c|c|c|c|}
\hline \lambda & (1,1,1,1,1,1) & (2,1,1,1,1)
&(2,2,1,1) & (2,2,2) & (3,1,1,1) & (3,2,1)\\
\hline (u_\lambda^{\cc})_0 & -\frac{1}{6} & -5 & -\frac{2693}{120} & -9 &
-\frac{1007}{60} & -\frac{907}{15} \\
\hline
\end{array}\\
&\begin{array}{|c||c|c|c|c|c|}
\hline\lambda & (3,3) & (4,1,1) & (4,2) & (5,1) & (6) \\
\hline (u_\lambda^{\cc})_0 & -\frac{50}{3} & -\frac{2129}{60} & -35 & -42 &
-22\\\hline
\end{array}
\end{align*}

\item Rank $3$:

\begin{displaymath}
\begin{array}{|c||c|c|c|c|c|c|}
\hline\lambda &(1) & (1,1) & (2) & (1,1,1) & (2,1) & (3) \\
\hline (u_\lambda^{\cc})_0 & 1 & -\frac{3}{2} & -\frac{3}{2}  & 4 & 10 & 4
\\\hline
\end{array}
\end{displaymath}

\begin{displaymath}
\begin{array}{|c||c|c|c|c|c|}
\hline \lambda & (1,1,1,1) & (2,1,1) & (2,2) & (3,1) & (4) \\
\hline (u_\lambda^{\cc})_0 & -\frac{111}{8} & -\frac{243}{4} & -\frac{75}{4}
&-42 & -\frac{55}{4}\\\hline
\end{array}
\end{displaymath}

\begin{displaymath}
\begin{array}{|c||c|c|c|c|c|c|c|}
\hline \lambda & (1,1,1,1,1) & (2,1,1,1) & (2,2,1) & (3,1,1) & (3,2) & (4,1) &
(5) \\
\hline (u_\lambda^{\cc})_0 & \frac{553}{10} & \frac{3553}{10} & \frac{6051}{20}
& \frac{693}{2} & 168 & 198 & \frac{273}{5}\\\hline
\end{array}
\end{displaymath}

\item Rank $4$:

\begin{displaymath}
\begin{array}{|c||c|c|c|c|c|c|}
\hline \lambda & (1) & (1,1) & (2) & (1,1,1) & (2,1) & (3)\\
\hline (u_\lambda^{\cc})_0 & 1 & -3 & -2 & 18 & 28 & \frac{22}{3} \\\hline
\end{array}
\end{displaymath}

\begin{displaymath}
\begin{array}{|c||c|c|c|c|c|}
\hline \lambda & (1,1,1,1) & (2,1,1) & (2,2) & (3,1) & (4) \\
\hline (u_\lambda^{\cc})_0 & -145 & -379 & -\frac{147}{2} & -165 & -35\\\hline
\end{array}
\end{displaymath}

\item Rank $5$:

\begin{displaymath}
\begin{array}{|c||c|c|c|c|c|c|}
\hline \lambda & (1) & (1,1) & (2) & (1,1,1) & (2,1) & (3)\\
\hline (u_\lambda^{\cc})_0 & 1 & -5 & -\frac{5}{2} & \frac{160}{3} & 60 &
\frac{35}{3}\\\hline
\end{array}
\end{displaymath}
\end{itemize}

\begin{remark}\text{}
\begin{itemize}
\item One notes that the signs are alternate: each partition $\lambda$ comes
with a sign $(-1)^{|\lambda|-1}$.

\item Looking at the partitions $(k)$, one recovers the coefficient
$\frac{(-1)^{k-1}}{k^2}\binom{r\cdot k}{k-1}$ obtained in Formula
(\ref{cherntauttrivial}).

\item The coefficient for a partition $(k,1)$ for $k\geq 2$, seems
 to be $\frac{(-1)^{k}(r-1)}{k+1}\binom{r\cdot k}{k}$.

\item In the rank $2$ case, it seems that
$u_{(1^k)}^{\cc}=\frac{(-1)^{k-1}}{k}(1+c_1+c_2)$ for $k\geq 1$.

\item For partitions of lengths two, the results are compatible with those of
\cite{NW}.
\end{itemize}
\end{remark}

\section{Other multiplicative characteristic classes}
\label{s:multTaut}

The result of Proposition \ref{cherntauttrivial} generalizes to all
multiplicative characteristic classes. By the splitting principle, any
multiplicative characteristic class is uniquely determined by its value on line
bundles, \ie by a power series $\phi(x)\in 1+x\IQ[[x]]$. Define from $\phi$ a
new power series $\psi(t)=\sum\limits_{k\geq 1}\psi_kt^k\in t\IQ[[t]]$ by the
relation:
\begin{displaymath}
\frac{\partial \psi}{\partial
t}\left(\frac{x}{\phi(-x)}\right)=\phi(-x).
\end{displaymath}
Then:

\begin{formula}[Boissi\`ere \& Nieper-Wi\ss kirchen {\cite[Theorem 4]{BNW}}]
Let $F$ be the trivial bundle of rank $r$ on $\IC^2$. Then:
\begin{displaymath}
\Phi(F)=\exp \left(\sum_{k\geq 1}
\frac{\psi_k}{k}\kq_{k}(1_{\IC^2})\right)\vac.
\end{displaymath}
\end{formula}

\part{Characteristic classes of the tangent bundle}
\label{p:tangent}

In this part, we study the characteristic classes of the tangent bundle on
$S^{[n]}$. We first recall the theoretic results. The study of the tangent
bundle is related to the study of the tautological bundles: The latter are
needed for example in the recursive computation of the Chern character (\S
\ref{ss:recursionChernT}).

\section{Shape of the formulas}

The first result concerns multiplicative characteristic classes:

\begin{theorem}[Boissi\`ere \cite{B}, Boissi\`ere \&
Nieper-Wi\ss kirchen \cite{BNW}]\label{th:mult} Let $\phi$ be a multiplicative
characteristic class. There are unique rational constants
$a_\lambda,b_\lambda,c_\lambda,d_\lambda$ such that generating series of the
$\phi$-classes of the tangent bundle on $S^{[n]}$ is given by:
\begin{displaymath}
\Phi(S):=\exp\left(\sum_{\lambda\in \cP} a_\lambda \kq_\lambda(1_S)+ b_\lambda
\kq_\lambda(e_S)+ c_\lambda \kq_\lambda(K_S)+ d_\lambda
\kq_\lambda(K_S^2)\right)\vac.
\end{displaymath}
\end{theorem}

The case of the Chern character is particular but the structure is similar:

\begin{theorem}[Boissi\`ere {\cite{B}}] There are unique rational constants
 $\alpha_\lambda,\beta_\lambda,\gamma_\lambda,\delta_\lambda$ such that the Chern character of the tangent bundle
 on $\Hilb(S)$ is given by:
\begin{displaymath}
\Ch(S):=\left(\sum_{\lambda\in \cP} \alpha_\lambda
\kq_\lambda(1_S)+\beta_\lambda \kq_\lambda(e_S)+ \gamma_\lambda
\kq_\lambda(K_S)+\delta_\lambda \kq_\lambda(K_S^2)\right)\unit.
\end{displaymath}
\end{theorem}

\section{The Chern class}

As a first step in the determination of the constants for multiplicative
characteristic classes, we consider the special case of the Chern class. We
wish to derive the constants in the universal formula:

\begin{equation}\label{gen:chernT}
\CC(S)=\exp\left(\sum_{\lambda\in \cP} a_\lambda \kq_\lambda(1_S)+ b_\lambda
\kq_\lambda(e_S)+ c_\lambda \kq_\lambda(K_S)+ d_\lambda
\kq_\lambda(K_S^2)\right)\vac.
\end{equation}

\subsection{Towards the $(1_S)$-series}\text{}

\begin{proposition}\label{prop:1SchernT} For $k\geq 0$ it is:
\begin{displaymath}
a_{2k+2}=0,\quad a_{2k+1}=\frac{(-1)^kC_k}{2k+1},
\end{displaymath}
where $C_k:=\frac{1}{k+1}\binom{2k}{k}$ is the $k$-th Catalan number.
\end{proposition}

\begin{proof} This result is proved in Boissi\`ere \cite[Theorem 1.1]{B}. We recall
briefly the main argument. Assume that the surface is the affine plane:
$S=\IC^2$. Then the formula for the Chern class is more simple: the canonical
class and the Euler class are zero, and all operators $\kq_\lambda$ for a
partition $\lambda$ of length $\ell(\lambda)>1$ are also zero. So the sum in
the exponential involves only the $1_S$-series, and in this series only the
partitions of length one. This gives all coefficients $a_\lambda$ for
$\lambda=(k)$, $k\geq 1$ since the Chern class takes the form:
\begin{displaymath}
\CC(\IC^2)=\exp\left(\sum_{k\geq 1} a_{k} \kq_{k}(1_S)\right)\vac.
\end{displaymath}
The computation is done by use of the equivariant cohomology of the Hilbert
scheme $(\IC^2)^{[n]}$ for the natural action of the torus $\IC^*$:

\begin{formula}[Boissi\`ere {\cite{B}}]\label{formula1} The Chern class is:
\begin{displaymath}
\CC(\IC^2)=\exp\left( \sum_{k\geq 0} \frac{(-1)^k C_k}{2k+1}
\kq_{2k+1}(1_S)\right)\vac,
\end{displaymath}
where $C_k:=\frac{1}{k+1}\binom{2k}{k}$ is the $k$-th Catalan number.
\end{formula}
This gives the constants as announced.
\end{proof}

\subsection{The $(1^k)$-series}\text{}

We give here all constants $a_{(1^k)},b_{(1^k)},c_{(1^k)},d_{(1^k)}$  ($k\geq
1$) in Formula (\ref{gen:chernT}) of the Chern class. Firstly, an application
of a result of G\"ottsche \cite{G} gives the series $b,d$:

\begin{proposition}
For $k\geq 1$ we have:
\begin{enumerate}
\item $b_{(1^k)}=\frac{1}{k}\sum_{i|k}i=:\frac{1}{k}\sigma_1(k)$;

\item $d_{(1^k)}=0$.
\end{enumerate}
\end{proposition}

\begin{proof}\text{}
\begin{enumerate}
\item We specialize the general formula to the case when $S$ is a K3 surface.
Since $K_S=0$, all terms involving the classes $K_S$ or $K_S^2$ disappear. The
total Chern class takes the form:
\begin{displaymath}
\CC(S)=\exp\left(\sum_{\lambda\in \cP} a_\lambda \kq_\lambda(1_S)+ b_\lambda
\kq_\lambda(e_S)\right)\vac.
\end{displaymath}

Since the cohomological degree of an operator $\kq_\lambda(\alpha)$ is:
\begin{displaymath}
\deg \kq_\lambda(\alpha) = 2(|\lambda|+\ell(\lambda))+|\alpha|-4,
\end{displaymath}
the only way to get in conformal weight $n$ a class of maximal degree $4n$ is
to use the operators $\kq_{(1^k)}(e_S)$. This means that:
\begin{displaymath}
\sum_{n\geq 0} e(S^{[n]})=\exp\left(\sum_{k\geq 1} b_{(1^k)}
\kq_{(1^k)}(e_S)\right)\vac.
\end{displaymath}
Denote by $\chi_S$ the Euler characteristic of $S$: $\chi_S=\int_Se_S$ or
equivalently $e_S=\chi_S x$ where $x$ denotes the cohomology class of a point.
Since
$$
\Delta^k_!e_S=\frac{1}{\chi_S^{k-1}}e_S\otimes\cdots\otimes e_S,
$$
we get $\kq_{(1^k)}(e_S)=\chi_S\kq_1(x)^k$ with
$\int_{S^{[k]}}\kq_1(x)^k\vac=1$. This implies that:
\begin{displaymath}
\sum_{n\geq 0} \int_{S^{[n]}}e(S^{[n]})t^n=\exp\left( \chi_S\sum_{k\geq 1}
b_{(1^k)}t^k\right).
\end{displaymath}
Now:
\begin{displaymath}
\sum_{n\geq 0} \int_{S^{[n]}}e(S^{[n]})t^n=\sum_{n\geq 0} \dim H^*(S^{[n]})t^n
\end{displaymath}
and using G\"ottsche's formula \cite{G} we get:
\begin{align*}
\sum_{n\geq 0} \int_{S^{[n]}}e(S^{[n]})t^n&=\prod_{m\geq 1} \left(\frac{1}{1-t^m}\right)^{\chi_S}\\
&=\exp\left(-\chi_S \sum_{m\geq 1}\ln(1-t^m)\right).
\end{align*}
This gives the relation $\sum\limits_{m\geq 1}\ln(1-t^m)=\sum\limits_{k\geq 1}
b_{(1^k)}t^k$ hence the result.

\item We make no assumption on the surface $S$, so the formula contains all
terms:
\begin{displaymath}
\CC(S)=\exp\left(\sum_{\lambda\in \cP} a_\lambda \kq_\lambda(1_S)+ b_\lambda
\kq_\lambda(e_S)+c_\lambda \kq_\lambda(K_S)+d_\lambda
\kq_\lambda(K_S^2)\right)\vac.
\end{displaymath}
Since the operators $\kq_{\lambda}(K_S^2)$ have same cohomological degree as
the operators $\kq_{\lambda}(e_S)$, the same argument as in the first assertion
gives in this case:
\begin{displaymath}
\sum_{n\geq 0} e(S^{[n]})=\exp\left(\sum_{k\geq 1} b_{(1^k)}
\kq_{(1^k)}(e_S)+d_{(1^k)} \kq_{(1^k)}(K_S^2)\right)\vac
\end{displaymath}
but since we have already obtained that
\begin{displaymath}
\sum_{n\geq 0} e(S^{[n]})=\exp\left(\sum_{k\geq 1} b_{(1^k)}
\kq_{(1^k)}(e_S)\right)\vac
\end{displaymath}
the operators $\kq_{(1^k)}(K_S^2)$ can not contribute to the Euler classes.
This forces the vanishing $d_{(1^k)}=0 $ for all $k\geq 1$.
\end{enumerate}
\end{proof}

In order to get the two series $a,c$, we make use of a result of Ohmoto
\cite{O} (in fact, this method recovers the series $b$ and $d$ obtained in the
preceding proposition, since the new argument uses a generalization of G{\"o}ttsche's formula).

\begin{proposition} For $k\geq 1$ it is:
\begin{displaymath}
a_{(1^k)}=-c_{(1^k)}=\frac{1}{k}\sum_{i|k}i=\frac{1}{k}\sigma_1(k).
\end{displaymath}
\end{proposition}

\begin{proof}
We follow Ohmoto \cite[Remark 2.4]{O}. Set $S^{(n)}:=S^n/\kS_n$, the quotient
of $S^n$ by the permutation action of the symmetric group $\kS_n$. The
\emph{Hilbert-Chow} morphism $\pi:S^{[n]}\rightarrow S^{(n)}$ is a crepant
resolution of singularities. The composite morphism
$S\xrightarrow{\Delta^n}S^n\rightarrow S^{(n)}$ induces a map
$\Delta^{(n)}:H^*(S)\rightarrow H^*(S^{(n)})$ (this is the Poincar\'e dual of
Ohmoto's $D^n$). By dimensional reasons, the cohomological push-forward
$f_!:H^*(S^{[n]})\rightarrow H^*(S^{(n)})$ vanishes on classes containing ar
least one operator $\kq_i$ with $i\geq 2$, and by definition
$f_!\kq_{(1^n)}(\alpha)\vac=\Delta^{(n)}\alpha$ for $\alpha\in H^*(S)$.

The total cohomology space $\bigoplus_{n\geq 0}H^*(S^{(n)})$ is equipped with a natural product
$\odot:~H^*(S^{(k)})\times~ H^*(S^{(l)})\rightarrow~ H^*(S^{(k+l)})$ (see \cite[\S 3.1]{O}) such that:
\begin{displaymath}
f_!\left(\kq_{(1^{k_1})}(\alpha_1)\cdots\kq_{(1^{k_r})}(\alpha_r)\vac\right)=\Delta^{(k_1)}(\alpha_1)\odot
\cdots\odot \Delta^{(k_r)}(\alpha_r).
\end{displaymath}
This gives us the image of the generating series of Formula (\ref{gen:chernT}):
\begin{displaymath}
f_!\CC(S)=\exp\left(\sum\limits_{k\geq 1}
a_{(1^k)}\Delta^{(k)}(1_S)+b_{(1^k)}\Delta^{(k)}(e_S)+c_{(1^k)}\Delta^{(k)}(K_S)+d_{(1^k)}\Delta^{(k)}(K^2_S)\right)
\end{displaymath}
where the exponential has to be taken for the $\odot$-ring structure. Now we apply Ohmoto's formula \cite[Formula (3)]{O}:
\begin{align*}
f_!\CC(S)&=\prod_{i\geq 1} (1-\Delta^{(i)})^{-\cc(S)}\\
&= \prod_{i\geq 1} \exp(-\ln(1-\Delta^{(i)})(\cc(S))\\
&= \exp \left(\sum_{i\geq 1} \sum_{j\geq 1} \frac{1}{j}\Delta^{(i\cdot j)}(1_S-K_S+e_S)\right)\\
&= \exp \left(\sum_{k\geq 1}
\sigma_1(k)\left(\Delta^{(k)}(1_S)+\Delta^{(k)}(e_S)-\Delta^{(k)}(K_S)\right)\right).
\end{align*}

This  gives $a_{(1^k)}=b_{(1^k)}=\sigma(k)$, $c_{(1^k)}=-\sigma(k)$ and $d_{(1^k)}=0$.
\end{proof}

\section{Other multiplicative classes}

Generalizing the Chern class to a multiplicative class $\phi$ gives other
series of universal coefficients:
\begin{displaymath}
\Phi(S)=\exp\left(\sum_{\lambda\in \cP} a_\lambda \kq_\lambda(1_S)+ b_\lambda
\kq_\lambda(e_S)+ c_\lambda \kq_\lambda(K_S)+ d_\lambda
\kq_\lambda(K_S^2)\right)\vac.
\end{displaymath}

\subsection{Towards the $(1_S)$-series}\text{}
\label{ss:1SHilb}

We proceed as in the Chern class case. For $S=\IC^2$, the $\phi$-class takes
the form:
\begin{displaymath}
\Phi(\IC^2)=\exp\left(\sum_{k\geq 1} a_{k} \kq_{k}(1_S)\right)\vac.
\end{displaymath}
This computation has been done in Boissi\`ere \& Nieper-Wi\ss kirchen
\cite{BNW}. By the splitting principle, any characteristic class is uniquely
determined by its value on a line bundle, \ie by a power series $\phi(x)\in
1+xA[[x]]$. Define from $\phi$ a new power series $\psi(t)=\sum\limits_{k\geq
1}\psi_kt^k\in tA[[t]]$ by the relation:
\begin{displaymath}
\frac{\partial \psi}{\partial
t}\left(\frac{x}{\phi(x)\phi(-x)}\right)=\phi(x)\phi(-x).
\end{displaymath}
Then:
\begin{formula}[Boissi\`ere \& Nieper-Wi\ss kirchen {\cite[Theorem 4]{BNW}}]
\begin{displaymath}
\Phi(\IC^2)=\exp\left( \sum_{k\geq 1} \frac{\psi_k}{k} \kq_{k}(1_S)\right)\vac.
\end{displaymath}
\end{formula}

\begin{remark}
The series $\psi$ is odd, so all coefficients $\psi_{2k}$ are zero.
\end{remark}

This general formula contains some nice special cases:

\noindent$\bullet$ For $\phi(x)=1+x$, one gets the Chern class (Formula
\ref{formula1}).

\noindent$\bullet$ For $\phi(x)=\frac{1}{1+x}$, one gets the Segr\'e class:
\begin{formula}[Boissi\`ere \& Nieper-Wi\ss kirchen {\cite[Example 6]{BNW}}]
\begin{displaymath}
\Se(\IC^2)=\exp\left( \sum_{k\geq 0} \frac{1}{(2k+1)^2}\binom{3k}{k}
\kq_{2k+1}(1_S)\right)\vac.
\end{displaymath}
\end{formula}

\noindent$\bullet$ For $\phi(x)=\sqrt{\frac{x}{1-\exp(-x)}}$ one gets the
square root of the Todd class:
\begin{formula}[Boissi\`ere \& Nieper-Wi\ss kirchen {\cite[Example 7]{BNW}}]
\begin{displaymath}
(\sqrt{\Td})(\IC^2)=\exp\left( \sum_{k\geq 0} \frac{1}{4^k\cdot (2k+1)\cdot
(2k+1)!} \kq_{2k+1}(1_S)\right)\vac.
\end{displaymath}
\end{formula}

\section{The Chern character}

We now consider the formula for the Chern character:

\begin{equation}\label{gen:cherncharT}
\Ch(S)=\left(\sum_{\lambda\in \cP} \alpha_\lambda
\kq_\lambda(1_S)+\beta_\lambda \kq_\lambda(e_S)+ \gamma_\lambda
\kq_\lambda(K_S)+\delta_\lambda \kq_\lambda(K_S^2)\right)\unit.
\end{equation}

\subsection{Towards the $(1_S)$-series}\text{}

\begin{proposition} For $k\geq 0$ it is:
\begin{displaymath}
\alpha_{2k+2}=0,\quad \alpha_{2k+1}=\frac{2}{(2k+1)!}.
\end{displaymath}
\end{proposition}

\begin{proof} This is proved in Boissi\`ere \cite[Theorem 1.1]{B}. The
argument is similar to the case of the Chern class (see Proposition
\ref{prop:1SchernT}) and the result is contained in the following formula:

\begin{formula}[Boissi\`ere {\cite[Theorem 1.1]{B}}]
\begin{displaymath}
\mathrm{Ch}(\IC^2)=\left(\sum_{k\geq 0}
\frac{2}{(2k+1)!}\kq_{2k+1}(1_S)\right)\unit.
\end{displaymath}
\end{formula}
This gives the constants as announced.
\end{proof}

\subsection{Towards the complete series}\text{}
\label{ss:recursionChernT}

Denote by $\kch T\in \End(\IH_S)$ the operator acting by multiplication by
$\ch(T^n_S)$ on each component of conformal weight $n$. In order to get
information on the series, we proceed to an implementation of the recursive
formula of Boissi\`ere \cite[Lemma 3.12]{B}\footnote{In the proof of
\cite[Proposition 3.10]{B}, the assumption $\int_S b_ib_j\td(S)=\delta_{i,j}$
should be only for cohomology with complex coefficients, otherwise one should
write $\int_Sb_ib_j\td(S)=\kappa_i\delta_{i,j}$ for some $\kappa_i\in \IQ$.
This does not affect the proof (just add the $\kappa_i$'s) since
$\Delta^2_!(\td(S))^{-1})=\ch(\cO_{\Xi^1})=\ch(\cO_{\Delta})$. Note that there
is an inaccuracy in the text since $\td(S)$ should be $\td(S)^{-1}$ at the end
of the proof.}$\text{}^{,}$\footnote{There is a typo in a computation on page
776: $\ch(\cO_S-T_S+\omega_S^\vee)=e_S$.}:
\begin{align*}
[\kch T,\kq_1(1_S)]=&\sum_\nu \frac{1}{\nu!} \kq_1^{(\nu)}(1_S) \\
&-\sum_{\nu} \frac{1}{\nu!} (\kq_1^{(\nu)}\circ \kG^\vee)\Delta^2_!(\td(S)^{-1}) \\
&+\sum_\nu \frac{(-1)^\nu}{\nu!} \kq_1^{(\nu)}(\exp(-K_S)) \\
&-\sum_{\nu} \frac{(-1)^\nu}{\nu!} (\kq_1^{(\nu)}\circ\kG)\Delta^2_!(\exp(-K_S)\td(S)^{-1}) \\
&-\kq_1(e_S)
\end{align*}
with $\td(S)^{-1}=1+\frac{K_S}{2}+\frac{2K_S^2-e_S}{12}$.

For the implementation, the computation with the $\kG$'s is explained in section \ref{chtautimplement}, and the case of the $\kG^\vee$'s is similar since it is easy to deduce from the results on $\kG$ the following commutation relation:
\begin{displaymath}
[\kG^\vee(\alpha),\kq_1(1_S)]=\exp(-\ad \kd)(\kq_1(\alpha)).
\end{displaymath}
This yields the following recursive formula:
\begin{align*}
\ch(T_n^S)=&\frac{1}{n}\kq_1(1_S)\ch(T_{n-1}^S)-\frac{1}{n!}\kq_1(e_S)\kq_1(1_S)^{n-1}\vac\\
&+\frac{1}{n!}\sum_{\nu=0}^{2n} \frac{1}{\nu!} \Bigg( \kq_1^{(\nu)}(1_S)\kq_1(1_S)^{n-1} \\
&\hspace{2cm}-(\kq_1^{(\nu)}\circ \kG^\vee)\Delta^2_!\left(1_S+\frac{K_S}{2}+\frac{2K_S^2-e_S}{12}\right)\kq_1(1_S)^{n-1}\\
&\hspace{2cm}+(-1)^\nu \kq_1^{(\nu)}\left(1_S-K_S+\frac{K_S^2}{2}\right)\kq_1(1)^{n-1} \\
&\hspace{2cm}- (-1)^\nu
(\kq_1^{(\nu)}\circ\kG)\Delta^2_!\left(1_S-\frac{K_S}{2}+\frac{2K_S^2-e_S}{12}\right)\kq_1(1_S)^{n-1}\Bigg)\vac.
\end{align*}

One gets the following series inside the brackets of Formula
(\ref{gen:cherncharT}):
\begin{displaymath}
\begin{array}{|c||c|c|c|c|c|c|}
\hline\lambda & (1) & (1,1) & (2) & (1,1,1) & (2,1) & (3) \\
\hline\alpha & 2 & -\frac{3}{2} & 0 & \frac{5}{9} & 0 & \frac{1}{3}\\
\hline\beta & -1 & \frac{5}{8} & 0 & -\frac{7}{27} & 0 & -\frac{5}{36}\\
\hline\gamma & -1 & -\frac{5}{12} & -1 & -\frac{1}{36} & \frac{13}{12} & -\frac{1}{6}\\
\hline\delta & \frac{1}{2} & \frac{1}{4} & \frac{1}{2} & \frac{53}{270} &
-\frac{5}{24} & \frac{4}{9}\\\hline
\end{array}
\end{displaymath}

\begin{remark}
\item The zeros are no surprise. In fact, for each partition $\lambda$ such
that $|\lambda|+\ell(\lambda)$ is odd, $\alpha_\lambda=\beta_\lambda=0$ since
if $S$ is a non-compact symplectic surface, these terms would contribute to the
even part of the Chern character, which is zero. The results are compatible
with those of \cite{NW}.
\end{remark}

\part{An implementation with {\sc Maude}}
\label{p:maude}

For documentation about {\sc Maude}, see \cite{Maude} or
\texttt{http://maude.cs.uiuc.edu/}.

\section{The script}

{\tiny
\begin{verbatim}
--- The total cohomology of Hilbert schemes is the module HILB

mod HILB is

    protecting RAT . --- Cohomology over the rational numbers

    sorts End Surf Part Alg . --- Endomorphims, Surface, Partitions, Algebra
    subsort Nat < Part .

    vars f g h k : End .
    vars a b : Rat .
    vars i j : Nat .
    vars m n : NzNat .
    vars p q : Part .
    vars c d s t u v : Surf .

-------------------
--- basic operators
-------------------

    --- Kronecker Symbol

    op kro(_ _) : Nat Nat -> Nat .

    eq kro (i j) = if i == j then 1 else 0 fi .

    --- Factorial

    op fact(_) : Nat -> Nat [memo].

    eq fact(0) = 1 .
    eq fact(i) = i * fact(i - 1) .

    --- Concatenation of partitions

    op _;_ : Part Part -> Part [assoc comm].

-----------------------------------------------------
--- Algebra of endomorphims over the rational numbers
-----------------------------------------------------

    ---- operations

    op O : -> End [ctor] . --- zero morphism
    op Id : -> End [ctor] . --- identity morphism

    op _+_ : End End -> End [assoc comm ctor prec 35 format (nt d d d)] . --- addition
    op __ : Rat End -> End [ctor prec 33] . --- multiplication by a scalar
    op _._ : End End -> End [ctor assoc prec 30] . --- composition
    op [__] : End End -> End . --- Lie bracket
    op P(__) : End Nat -> End . --- Iterated composition (power)

    --- Axioms for the addition

    eq f + f = 2 f .
    eq O + f = f .

   --- Axioms for the multiplication by a scalar

    eq a f + b f = (a + b) f .
    eq a (f + g) = a f + a g .
    eq a f + f = (a + 1) f .
    eq a (b f) = (a * b) f .
    eq 1 f = f .
    eq 0 f = O .
    eq a O = O .

    --- Axioms for the composition

    eq (a f) . g = a (f . g) .
    eq f . (a g) = a (f . g) .
    eq f . (h + k) = f . h + f . k .
    eq (f + g ) . h = f . h + g . h .
    eq Id . f = f .
    eq f . Id = f .
    eq O . f = O .
    eq f . O = O .

    --- Axioms for the Lie bracket

    eq [ f g ] = f . g + (-1) g . f .

    --- Axiom for iterated composition (power)

    eq P(f 0) = Id .
    eq P(f i) = f . P( f (i - 1) ) .


--------------------------------------------
--- Rational cohomology algebra of a surface
--------------------------------------------

   --- Usual classes

   op o : -> Surf [ctor] . --- zero class
   op I : -> Surf [ctor] . --- unit class
   op K : -> Surf [ctor] . --- canonical class
   op e : -> Surf [ctor] . --- Euler class
   op C : -> Surf . --- a degree 2 class (first Chern class)
   op D : -> Surf . --- a degree 4 class (second Chern class)

   --- Operations

   op _+_ : Surf Surf -> Surf [assoc comm ctor prec 35] . --- addition
   op __ : Rat Surf -> Surf [ctor prec 33] . --- external law
   op _._ : Surf Surf -> Surf [ctor assoc comm prec 30] . --- cup product
   op int(__) : Surf Surf -> Rat [comm] . --- intersection product
   op [_'] : Surf -> Surf . --- simulation of Kuenneth decomposition

   --- Axioms for the addition

    eq s + s = 2 s .
    eq o + s = s .

   --- Axioms for the external law

    eq a s + b s = (a + b) s .
    eq a (s + t) = a s + a t .
    eq a s + s = (a + 1) s .
    eq a (b s) = (a * b) s .
    eq 1 s = s .
    eq 0 s = o .
    eq a o = o .

    --- Axioms for the cup product

    eq (a s) . t = a (s . t) .
    eq s . (a t) = a (s . t) .
    eq s . (t + u) = s . t + s . u .
    eq (s + t) . u = s . u + t . u .
    eq I . s = s .
    eq o . s = o .

   --- special cases (for K3 or abelian surface)

    --- eq K = o .
    --- eq e = o .

   --- universal relations

   eq K . K . K = o .
   eq e . e = o .
   eq e . K = o .
   eq C . C . C = o .
   eq C . C . K = o .
   eq C . K . K = o .
   eq C . e = o .
   eq C . D = o .
   eq D . D = o .
   eq D . K = o .
   eq D . e = o .

--------------------
--- Vertex operators
--------------------

   op |> : -> End [ctor]. --- the vaccum, considered as an operator for simplicity
   op d : -> End . --- boundary operator
   op q(___) : Nat Surf NzNat -> End . --- derived Nakajima operator
   op <__> : Surf NzNat -> End . --- non-derived Nakajima operator
   op <<__>> : Surf Part -> End . --- Nakajima operator for a partition
   op L(__) : Surf NzNat -> End . --- Virasoro operator
   op Rec(___) : Surf NzNat NzNat -> End . --- recursive computation of L(__).|>
   op CH(__) : Surf Nat -> End [memo]. --- Chern character of a tautological class
   op Rec2(___) : Surf Nat Nat -> End . --- recursive computation of CH(__).|>
   op CHv(__) : Surf Nat -> End [memo]. --- Chern character of a dual tautological class
   op Rec2v(___) : Surf Nat Nat -> End . --- recursive computation of CHv(__).|>
   op C1(__) : Surf Nat -> End [memo]. --- Chern class of a tautological rank 1 bundle
   op C2(___) : Surf Surf Nat -> End [memo]. --- Chern class of a tautological rank 2 bundle
   op C2t(_) : Nat -> End [memo]. --- Chern class of a tautological rank 2 trivial bundle
   op C3t(_) : Nat -> End [memo]. --- Chern class of a tautological rank 3 trivial bundle
   op C4t(_) : Nat -> End [memo]. --- Chern class of a tautological rank 4 trivial bundle
   op C5t(_) : Nat -> End [memo]. --- Chern class of a tautological rank 5 trivial bundle
   op CHT(_) : Nat -> End [memo]. --- Chern character of the tangent bundle
   op Rec3(___) : Nat Nat NzNat -> End . --- recursive computation of CHT(_).|>
   --- Basics

   eq d . |> = O .
   eq < o n > = O .
   eq << o p >> = O .

   --- Derived Nakajima operators

   eq q(0 u n) = < u n > .
   eq q(i u n) = [ d  q( (i - 1) u n ) ] .
   eq d . < u n > = < u n > . d + n L( u n ) + ( n * (n - 1)) / 2 < (K . u) n > .

   --- Virasoro operators

   eq L(u n) . < v m > = < v m > . L(u n) + (-1) * m < (u . v) (n + m) > .
   eq L(u n) . |> = Rec(u n 1) .
   eq Rec(u n n) = O .
   eq Rec(u n m) = (1 / 2) < [ u '] m > . < [ u '] (n - m) > . |> + Rec(u n (m + 1)) .

   --- Chern character of a tautological class

   eq CH(c 0) = O .
   eq CH(c n) = (1 / n) < I 1 > . CH( c (n - 1) ) + (1 / fact(n) ) Rec2(c 0 ( (2 * n) ) ) . P( < I 1 > (n - 1) ) . |> .

   eq Rec2(c j j) = (1 / fact(j) ) q(j c 1) .
   eq Rec2(c i j) = (1 / fact(i) ) q(i c 1) + Rec2(c (i + 1) j) .

   --- Chern character of a dual tautological class

   eq CHv(c 0) = O .
   eq CHv(c n) = (1 / n) < I 1 > . CHv( c (n - 1) ) + (1 / fact(n) ) Rec2v(c 0 ( (2 * n) ) ) . P( < I 1 > (n - 1) ) . |> .

   eq Rec2v(c j j) = ( ( (-1) ^ j ) / fact(j) ) q(j c 1) .
   eq Rec2v(c i j) = ( ( (-1) ^ i ) / fact(i) ) q(i c 1) + Rec2v(c (i + 1) j) .

   --- Chern class of a tautological bundle

       --- Rank 1: c is the first Chern class of a line bundle
   eq C1(c 0) = |> .
   eq C1(c i) = (1 / i) (< c 1 > + q(1 I 1)) . C1(c (i - 1)) .

       --- Rank 2: c is the first Chern class, d the second Chern class
   eq C2(c d 0) = |> .
   eq C2(c d i) = (1 / i) (< I 1 > + < c 1 > + < d 1 > + 2 q(1 I 1)
                  + q(1 c 1) + q(2 I 1) ) . C2(c d (i - 1)) .

       --- special cases: trivial bundles of ranks 2,3,4,5
   eq C2t(0) = |> .
   eq C2t(i) = (1 / i) (< I 1 > + 2 q(1 I 1) + q(2 I 1) ) . C2t((i - 1)) .

   eq C3t(0) = |> .
   eq C3t(i) = (1 / i) ( < I 1 > + 3 q(1 I 1) + 3 q(2 I 1) + q(3 I 1) ) . C3t((i - 1)) .

   eq C4t(0) = |> .
   eq C4t(i) = (1 / i) ( < I 1 > + 4 q(1 I 1) + 6 q(2 I 1) + 4 q(3 I 1)
               + q(4 I 1)) . C4t((i - 1)) .

   eq C5t(0) = |> .
   eq C5t(i) = (1 / i) ( < I 1 > + 5 q(1 I 1) + 10 q(2 I 1) + 10 q(3 I 1)
               + 5 q(4 I 1) + q(5 I 1)) . C5t((i - 1)) .

   --- Chern character of the tangent bundle

   eq CHT(0) = O .
   eq CHT(n) = (1 / n) < I 1 > . CHT((n - 1))
               + ((-1) / fact(n)) < e 1 > . P( < I 1 > (n - 1) ) . |>
               + (1 / fact(n) ) Rec3( 0 (2 * n) n) .

   eq Rec3(j j n) = (1 / fact(j)) ( q(j I 1) . P( < I 1 > (n - 1) ) . |>
                    + (-1) fact(n - 1) q(j [ I '] 1) . CHv( [ I '] (n - 1) )
                    + ((-1) / 2) fact(n - 1) q(j [ K '] 1) . CHv( [ K '] (n - 1))
                    + ((-1) / 6) fact (n - 1) q(j [ (K . K) '] 1) . CHv([ (K . K) '] (n - 1) )
                    + (1 / 12) fact(n - 1) q(j [ e '] 1) . CHv( [ e '] (n - 1) )
                    + ((-1) ^ j) q(j I 1) . P( < I 1 > (n - 1) ) . |>
                    + ((-1) ^ (j + 1)) q(j K 1) . P( < I 1 > (n - 1) ) . |>
                    + (((-1) ^ j) / 2) q(j (K . K) 1) . P( < I 1 > (n - 1) ) . |>
                    + ((-1) ^ (j + 1)) fact(n - 1) q(j [ I '] 1) . CH( [ I '] (n - 1))
                    + (((-1) ^ j) / 2) fact(n - 1) q(j [ K '] 1) . CH( [ K '] (n - 1))
                    + (((-1) ^ (j + 1)) / 6) fact(n - 1) q(j [ (K . K) '] 1) . CH( [ (K . K) '] (n - 1))
                    + (((-1) ^ j) / 12) fact(n - 1) q(j [ e '] 1) . CH( [ e '] (n - 1)) ) .
   eq Rec3(i j n) = (1 / fact(i)) ( q(i I 1) . P( < I 1 > (n - 1) ) . |>
                    + (-1) fact(n - 1) q(i [ I '] 1) . CHv( [ I '] (n - 1))
                    + ((-1) / 2) fact(n - 1) q(i [ K '] 1) . CHv( [ K '] (n - 1))
                    + ((-1) / 6) fact(n - 1) q(i [ (K . K) '] 1) . CHv( [ (K . K) '](n - 1))
                    + (1 / 12) fact(n - 1) q(i [ e '] 1) . CHv( [ e '] (n - 1))
                    + ((-1) ^ i) q(i I 1) . P( < I 1 > (n - 1) ) . |>
                    + ((-1) ^ (i + 1)) q(i K 1) . P( < I 1 > (n - 1) ) . |>
                    + (((-1) ^ i) / 2) q(i (K . K) 1) . P( < I 1 > (n - 1) ) . |>
                    + ((-1) ^ (i + 1)) fact(n - 1) q(i [ I '] 1) . CH( [ I '] (n - 1))
                    + (((-1) ^ i) / 2) fact(n - 1) q(i [ K '] 1) . CH( [ K '] (n - 1))
                    + (((-1) ^ (i + 1)) / 6) fact(n - 1) q(i [ (K . K) '] 1) . CH( [ (K . K) '] (n - 1))
                    + (((-1) ^ i) / 12) fact(n - 1) q(i [ e '] 1) . CH( [ e '] (n - 1)) )
                    + Rec3( (i + 1) j n) .

   --- Rules for the final simplifications of Kuenneth decompositions

   rl [1] : < [ u ']  n > . < [ u '] m > => << u (n ; m) >> .
   rl [2] : < ( s . [ u '] )  n > . < [ u '] m > => << ( s . u ) (n ; m) >> .
   rl [3] : < [ u ']  n > . < ( t . [ u '] ) m > => << ( t . u ) (n ; m) >> .
   rl [4] : < ( s . [ u '] )  n > . < ( t . [ u '] ) m > => << ( s . t . u ) (n ; m) >> .

   rl [5] : <<  [ u ']  p >> . << [ u ']  q >> => <<  u  (p ; q) >> .
   rl [6] : << ( s . [ u '] )  p >> . << [ u ']  q >> => << ( s . u ) (p ; q) >> .
   rl [7] : << [ u ']  p >> . << ( t . [ u '] ) q >> => << ( t . u ) (p ; q) >> .
   rl [8] : << ( s . [ u '] )  p >> . << ( t . [ u '] ) q >> => << ( s . t . u ) (p ; q) >> .

   rl [9] : <<  [ u ']  p >> . < [ u ']  m > => <<  u  (p ; m) >> .
   rl [10] : << ( s . [ u '] )  p >> . < [ u ']  m > => << ( s . u ) (p ; m) >> .
   rl [11] : << [ u ']  p >> . < ( t . [ u '] ) m > => << ( t . u ) (p ; m) >> .
   rl [12] : << ( s . [ u '] )  p >> . < ( t . [ u '] ) m > => << ( s . t . u ) (p ; m) >> .

   rl [13] : < [ u ']   n > . << [ u ']  q >> => <<  u  (n ; q) >> .
   rl [14] : < ( s . [ u '] )  n > . << [ u ']  q >> => << ( s . u ) (n ; q) >> .
   rl [15] : < [ u ']  n > . << ( t . [ u '] ) q >> => << ( t . u ) (n ; q) >> .
   rl [16] : < ( s . [ u '] )  n > . << ( t . [ u '] ) q >> => << ( s . t . u ) (n ; q) >> .

   rl [17] : < ([u '] . [u ']) n > => < (e . u) n > .
   rl [18] : << ([u '] . [u ']) p >> => << (e . u) p >> .

   rl [19] : < (s . [u '] . [u ']) n > => < (s . e . u) n > .
   rl [20] : << (s . [u '] . [u ']) p >> => << (s . e . u) p >> .


   rl [1b] : <  [ u ']  n > . f . < [ u ']  m > =>  f . <<  u  (n ; m) >> .
   rl [2b] : < ( s . [ u '] )  n > . f . < [ u ']  m > => f . << ( s . u ) (n ; m) >> .
   rl [3b] : < [ u ']  n > . f . < ( t . [ u '] ) m > => f . << ( t . u ) (n ; m) >> .
   rl [4b] : < ( s . [ u '] )  n > . f . < ( t . [ u '] ) m > => f . << ( s . t . u ) (n ; m) >> .

   rl [5b] : <<  [ u ']  p >> . f . << [ u ']  q >> => f . <<  u  (p ; q) >> .
   rl [6b] : << ( s . [ u '] )  p >> . f . << [ u ']  q >> => f . << ( s . u ) (p ; q) >> .
   rl [7b] : << [ u ']   p >> . f . << ( t . [ u '] ) q >> => f . << ( t . u ) (p ; q) >> .
   rl [8b] : << ( s . [ u '] )  p >> . f . << ( t . [ u '] ) q >> => f . << ( s . t . u ) (p ; q) >> .

   rl [9b] : << [ u ']   p >> . f . < [ u ']  m > => f . <<  u  (p ; m) >> .
   rl [10b] : << ( s . [ u '] )  p >> . f . < [ u ']  m > => f . << ( s . u ) (p ; m) >> .
   rl [11b] : << [ u ']  p >> . f . < ( t . [ u '] ) m > => f . << ( t . u ) (p ; m) >> .
   rl [12b] : << ( s . [ u '] )  p >> . f . < ( t . [ u '] ) m > => f . << ( s . t . u ) (p ; m) >> .

   rl [13b] : <  [ u ']   n > . f . << [ u ']  q >> => f . <<  u  (n ; q) >> .
   rl [14b] : < ( s . [ u '] )  n > . f . << [ u ']  q >> => f . << ( s . u ) (n ; q) >> .
   rl [15b] : < [ u ']  n > . f . << ( t . [ u '] ) q >> => f . << ( t . u ) (n ; q) >> .
   rl [16b] : < ( s . [ u '] )  n > . f . << ( t . [ u '] ) q >> => f . << ( s . t . u ) (n ; q) >> .

endm

------------
--- Examples
------------

   --- Chern class of a tautological rank 2 bundle
   rew C2(c d 1) .
   rew C2(c d 2) .
   rew C2(c d 3) .

   --- Chern character of the tangent bundle
   rew CHT(1) .
   rew CHT(2) .
   rew CHT(3) .


\end{verbatim}}

\section{Final remarks}

Although we implemented lots of classes, things remain to do, in particular:
\begin{itemize}
\item Write and implement recursions for the Segre classes of tautological
bundles.

\item Write and implement a recursion for the Chern class of the tangent
bundle.

\item Find satisfactory models for the general terms of the series obtained by
symbolic computations. These models are still missing.
\end{itemize}

We decided not to implement further operators here. The methods are more or
less straightforward generalization of Lehn's ideas in \cite{Lehn}, but the
formulas would be very long: the complexity occurs in the decomposition of the
characteristic class of the tensor product of a vector bundle with a line
bundle.

\nocite{*}
\bibliographystyle{amsplain}
\bibliography{HilbClass}

\end{document}